\documentclass[12pt,leqno,draft]{article}

\def\diam{\mathop{\rm diam}}
\def\dist{\mathop{\rm dist}}

\newtheorem{theorem}{Theorem}
\newtheorem{lemma}[theorem]{Lemma}
\newtheorem{proposition}[theorem]{Proposition}
\newtheorem{definition}[theorem]{Definition}
\newtheorem{corollary}[theorem]{Corollary}

\newcommand{\begintheorem}{\addtocounter{equation}{1}\begin{theorem}}
\newcommand{\beginlemma}{\addtocounter{equation}{1}\begin{lemma}}
\newcommand{\beginproposition}{\addtocounter{equation}{1}\begin{proposition}}
\newcommand{\begindefinition}{\addtocounter{equation}{1}\begin{definition}}
\newcommand{\begincorollary}{\addtocounter{equation}{1}\begin{corollary}}

\begin{document}

\title{Cellular structures, quasisymmetric mappings, and spaces of
homogeneous type}

\author{Stephen Semmes \\
        Rice University}

\date{}

\maketitle

        Let $X$ be a compact Hausdorff topological space.  A
collection $\mathcal{C}$ of nonempty subsets of $X$ may be described
as a \emph{cellular structure} for $X$ if it satisfies the following
three properties.  First, each $C \in \mathcal{C}$ is both open and
closed, and $X \in \mathcal{C}$.  Second, $\mathcal{C}$ is a base for
the topology of $X$.  Third, if $C, C' \in \mathcal{C}$, then either
\begin{equation}
        C \cap C' = \emptyset, \hbox{ or } C \subseteq C',
                               \hbox{ or } C' \subseteq C.
\end{equation}
In this case, a set $C \in \mathcal{C}$ may be called a \emph{cell} in
$X$.  A compact Hausdorff space with a cellular structure is a
\emph{cellular space}.

        For example, the usual construction of the Cantor set leads to
a natural cellular structure, where the cells are the parts of the
Cantor set in the closed intervals generated in the construction.

        Of course, a Hausdorff topological space with a base for its
topology consisting of sets that are both open and closed is
automatically totally disconnected, in the sense that there are no
connected subsets with more than one element.  The collection of all
subsets of the space that are both open and closed is then an algebra
of sets as well as a base for the topology.  One can think of a
cellular structure as a kind of geometric structure on such a space.

        For the sake of simplicity, let us restrict our attention to
compact spaces, although one could consider non-compact spaces too.
For instance, one might consider locally compact Hausdorff spaces that
are $\sigma$-compact.

        Let $(X, \mathcal{C})$ be a cellular space, and suppose that
$A \subseteq X$ is open and closed.  In particular, $A$ is compact,
since $X$ is compact.  Because $A$ is open and $\mathcal{C}$ is a base
for the topology of $X$, $A$ can be expressed as the union of a
collection of cells.  If $A = \emptyset$, then one can interpret this
as meaning that $A$ is the union of the empty collection of cells.
This is an open covering of $A$, since cells are open sets.  By
compactness, $A$ is the union of finitely many cells.  Using the
nesting property of cells, it follows that $A$ is the union of
finitely many disjoint cells.

        Suppose that $C_1, \ldots, C_n$ are finitely many
pairwise-disjoint cells in $X$.  In particular, $C_1, \ldots C_n$ are
both open and closed, as is
\begin{equation}
        A = X \backslash (C_1 \cup \cdots \cup C_n).
\end{equation}
The preceding observation implies that $X$ is the union of a
collection of finitely many pairwise-disjoint cells that includes the
$C_i$'s.

        Suppose that $C_1, \ldots, C_n$ are finitely many
pairwise-disjoint cells whose union is $X$.  If $C$ is a cell such
that $C_i \subseteq C$ for some $i$, $1 \le i \le n$, then $C$ is the
union of some of the $C_j$'s.  Hence there can only be finitely many
such cells $C$.  It follows that every cell in $X$ is contained in
only finitely many other cells.

        If $C$ is a cell with at least two elements, then every point
in $C$ is contained in a cell that is a proper subset of $C$.  By
compactness, $C$ is the union of finitely many cells that are proper
subsets of $C$.  As usual, the smaller cells can also be taken to be
pairwise disjoint.  If a cell has only one element $p$, then $p$ is an
isolated point in $X$.

        For each $p \in X$,
\begin{equation}
        \mathcal{C}(p) = \{C \in \mathcal{C} : p \in C\}
\end{equation}
is linearly ordered by inclusion, because of the nesting property for cells.
Thus the collection of cells that contains a fixed cell $C_0$ is
finite and linearly ordered by inclusion, which means that any cell in
some collection of cells is contained in a maximal cell in the same
collection.  Using this, one can check that every cell $C$ with at
least two elements is the union of finitely many cells that are proper
subsets of $C$ and maximal with respect to inclusion.  Note that
maximal cells in any collection are automatically pairwise disjoint.

        A consequence of these remarks is that there are only finitely
or countably many cells in $X$.  This is trivial when $X$ has only one
element, and otherwise $X$ is the union of finitely many
pairwise-disjoint proper sub-cells $C_1, \ldots, C_n$.  There are only
finitely many cells that contain one of the $C_i$'s, and every other
cell is contained in one of the $C_i$'s.  Each $C_i$ with at least two
elements is also a union of finitely many pairwise-disjoint proper
sub-cells, and so one can repeat the process.  Every cell is contained
in only finitely many other cells, and hence is reached in finitely
many steps.

        Let $(X, \mathcal{C})$ be a cellular space, and suppose that
$Y \subseteq X$ is nonempty and compact.  Put
\begin{equation}
 \mathcal{C}_Y = \{C \cap Y : C \in \mathcal{C}, \, C \cap Y \ne \emptyset\}.
\end{equation}
This is a cellular structure on $Y$, which is induced from the one on
$X$.

        The example of the Cantor set can be extended, as follows.
Let $X_1, X_2, \ldots$ be a sequence of finite sets with at least two
elements, and let $X$ be the set of sequences $x = \{x_i\}_{i =
1}^\infty$ such that $x_i \in X_i$ for each $i$.  Thus $X$ is the
Cartesian product of the $X_i$'s, which is a compact Hausdorff space
with respect to the product topology using the discrete topology on
each $X_i$.  For every nonnegative integer $l$ and $x \in X$, let
$N_l(x)$ be the set of $y \in X$ such that $y_i = x_i$ when $i \le l$.
Note that $N_0(x) = X$, and $N_l(x)$ is both open and closed in $X$
for all $l \ge 0$ and $x \in X$.  The collection of $N_l(x)$'s is a
base for the product topology on $X$.  It is easy to check that the
collection of $N_l(x)$'s defines a cellular structure for $X$.  Let us
call this the \emph{product cellular structure} on $X$.

        If $(X, \mathcal{C})$ is a cellular space, then there is a
natural graph $\mathcal{T}$ whose vertices are the cells in $X$.
Specifically, we can attach an edge between the vertices associated to
two cells $C$, $C'$ when $C' \subseteq C$, $C' \ne C$, and $C'$ is a
maximal proper sub-cell in $C$.  The nesting property of cells implies
that $\mathcal{T}$ is a tree.  By definition, $X$ is a cell, which we
can take to be the root of the tree.

        Conversely, suppose that $\mathcal{T}$ is a locally-finite
tree with root $\tau$.  A \emph{ray} in $\mathcal{T}$ is a simple path
beginning at $\tau$ and continuing as long as possible.  More
precisely, a ray in $\mathcal{T}$ may stop after finitely many steps
when it arrives at a vertex with no additional edge to follow, or it
may traverse infinitely many edges.  Let $X$ be the set of rays in
$\mathcal{T}$.  It is convenient to represent a ray in $\mathcal{T}$
as an infinite sequence of vertices in $\mathcal{T}$, where the last
vertex of a finite ray is repeated indefinitely.  For each nonnegative
integer $l$, let $X_l$ be the set of vertices of $\mathcal{T}$ that
can be reached from $\tau$ in $\le l$ steps.  Thus $X_l$ has only
finitely many elements for every $l \ge 0$, and $X$ can be identified
with a subset of the Cartesian product of the $X_l$'s.  It is easy to
see that this is a closed set in the product topology, so that $X$
becomes a compact Hausdorff space using the induced topology.  Each
finite simple path in $\mathcal{T}$ starting at $\tau$ determines a
set of rays, i.e., the set of rays that are continuations of the path.
One can use these sets of rays as cells in $X$, which are the same as
the ones induced by the product space.

        Suppose that $(X, d(x, y))$ is a compact ultrametric space.
This means that $(X, d(x, y))$ is a compact metric space, and that
\begin{equation}
        d(x, z) \le \max(d(x, y), d(y, z))
\end{equation}
for every $x, y, z \in X$.  In an ultrametric space, closed balls with
positive radii are open sets, and form a base for the topology.  The
ultrametric version of the triangle inequality also implies that any
two closed balls are either disjoint or one is contained in the other,
so that closed balls with positive radii determine a cellular
structure on $X$.  More precisely, a set $C \subseteq X$ would be a
cell if it could be expressed as a closed ball with some center $x \in
X$ and radius $r > 0$, but $x$ and $r$ are not necessarily uniquely
determined by $C$.

        If $X = \prod_{i = 1}^\infty X_i$ with the product cellular
structure, then compatible ultrametrics on $X$ can be obtained as
follows.  Let $\rho = \{\rho_i\}_{i = 0}^\infty$ be a strictly
decreasing sequence of positive real numbers such that $\rho_0 = 1$ and
\begin{equation}
        \lim_{i \to \infty} \rho_i = 0.
\end{equation}
For each $x, y \in X$, put $d_\rho(x, y) = 0$ when $x = y$, and
otherwise
\begin{equation}
        d_\rho(x, y) = \rho_n
\end{equation}
where $l$ is the largest nonnegative integer such that $x_i = y_i$
when $i \le n$.  It is not difficult to verify that $d_\rho(x, y)$ is
an ultrametric on $X$ for which the corresponding topology is the
product topology, and for which the closed balls are the cells in the
product cellular structure.

        A standard regularity condition for $\rho$ asks that there
be real numbers $0 < a \le b < 1$ such that
\begin{equation}
\label{regularity condition}
        a \le \frac{\rho_{i + 1}}{\rho_i} \le b
\end{equation}
for each $i \ge 0$.  For instance, this holds if $\rho_i$ is the $i$th
power of a fixed positive real number less than $1$.  If $\rho$,
$\widetilde{\rho}$ are two such sequences, then the corresponding
ultrametrics $d_\rho(x, y)$, $d_{\widetilde{\rho}}(x, y)$ are
quasisymmetrically equivalent, in the sense that the identity mapping
on $X$ is quasisymmetric as a mapping from $(X, d_\rho(x, y))$ to $(X,
d_{\widetilde{\rho}}(x, y))$.

        Remember that a metric space $(X, d(x, y))$ is said to be
\emph{doubling} if every ball in $X$ can be covered by a bounded
number of balls of half the radius.  A positive Borel measure on $X$
is said to be a \emph{doubling measure} with respect to the metric if
the measure of every ball is bounded by a constant times the measure
of the ball with the same center and half the radius.  A well known
covering argument implies that a metric space with a doubling measure
is doubling.

        In analogy with this, let us say that a cellular space $(X,
\mathcal{C})$ is \emph{doubling} if there is a $k_1 \ge 1$ such that
every cell $C$ in $X$ contains no more than $k_1$ maximal proper
sub-cells.  Similarly, a positive Borel measure $\mu$ on $X$ is a
\emph{doubling measure} with respect to the cellular structure if
\begin{equation}
        0 < \mu(C) < \infty
\end{equation}
for every cell $C$, and if there is a $k_2 \ge 1$ such that
\begin{equation}
        \mu(C) \le \mu(C')
\end{equation}
whenever $C$, $C'$ are cells such that $C' \subseteq C$, $C' \ne C$,
and $C'$ is a maximal proper sub-cell in $C$.  This implies that $C$
has at most $k_2$ maximal proper sub-cells, since $C$ is the disjoint
union of its maximal proper sub-cells.

        For example, if $X = \prod_{i = 1}^\infty X_i$ with the
product cellular structure, then $X$ is doubling with constant $k_1$
if and only if each $X_i$ has at most $k_1$ elements.  A nice class of
measures on $X$ is given by product measures $\mu = \prod_{i =
1}^\infty \mu_i$, where each $\mu_i$ is a probability measure on the
finite set $X_i$.  Thus $\mu_i$ is defined by assigning weights to the
elements of $X_i$ whose sum is $1$, and $\mu$ is doubling on $X$ with
constant $k_2$ if and only if the $\mu_i$ measure of each element of
$X_i$ is at least $1/k_2$.

        If $\rho$ is a strictly decreasing sequence of positive real
numbers that satisfies the regularity condition (\ref{regularity
condition}), then the doubling condition for $X$ as a metric space
with the metric $d_\rho(x, y)$ is also equivalent to the boundedness
of the number elements of the $X_i$'s.  Also, a positive Borel measure
on $X$ is then doubling with respect to the metric $d_\rho(x, y)$ if
and only if it is doubling with respect to the associated cellular
structure.

        Let $(X, \mathcal{C})$ be a cellular space with a metric $d(x,
y)$ that determines the same topology on $X$.  A more abstract version
of the regularity condition (\ref{regularity condition}) for the
compatibility of the metric $d(x, y)$ with the cellular structure
$\mathcal{C}$ asks that there be positive real numbers $\alpha$,
$\beta$, $\gamma$ with $\alpha \le \beta < 1$ such that
\begin{equation}
\label{regularity condition, 1}
        \alpha \, \diam C \le \diam C' \le \beta \diam C
\end{equation}
for every cell $C$ in $X$ and maximal proper sub-cell $C'$ in $C$, and
\begin{equation}
\label{regularity condition, 2}
        \dist (C', C'') \ge \gamma \, \diam C
\end{equation}
when $C'$, $C''$ are distinct maximal proper sub-cells of a cell $C$.
As usual, $\diam C$ denotes the diameter of $C$, which is the supremum
of the distances between elements of $C$, and $\dist (C', C'')$
denotes the distance between $C'$ and $C''$, which is to say the
infimum of the distances between elements of $C'$ and $C''$.  In the
special case where $X = \prod_{i = 1}^\infty X_i$ with the product
cellular structure and $d(x, y) = d_\rho(x, y)$, $\alpha$ and $\beta$
correspond exactly to $a$ and $b$, and one can take $\gamma = 1$.

        Note that a cell $C$ with at least two elements has a proper
sub-cell, since the cells form a base for the topology.  Conversely,
if a cell $C$ has a proper sub-cell, then $C$ has at least two
elements, and hence the diameter of $C$ is positive.  If $C'$ is a
maximal proper sub-cell of a cell $C$, then the preceding regularity
condition implies that $C'$ has positive diameter as well.  Applying
this repeatedly, it follows that $X$ has no isolated points when $X$
has at least two elements.

        Suppose that $d(x, y)$ is an ultrametric on $X$ and
$\mathcal{C}$ consists of the closed balls in $X$ with positive
radius.  If $C$ is a cell in $X$ with diameter $r$, then $C$ is the
same as the closed ball in $X$ defined by $d(x, y)$ with radius $r$
and centered at any element of $C$.  It may be that $r = 0$, so that
$C$ consists of a single point $p$, in which case $p$ should be an
isolated point in $X$.  In any case, it may be possible to represent
$C$ as a ball of radius larger than $r$.  Note that the diameter of a
ball of radius $t \ge 0$ in an ultrametric space is less than or equal
to $t$, while in an ordinary metric space it is less than or equal to
$2 \, t$ and often equal to $2 \, t$.

        If $C$, $C'$ are cells in $X$ such that $C' \subseteq C$ and
$C' \ne C$, then
\begin{equation}
        \diam C' < \diam C.
\end{equation}
Indeed, $\diam C' \le \diam C$ since $C' \subseteq C$, and equality of
the diameters would imply that $C' = C$, by the previous remarks.  If
$C'$, $C''$ are cells in $X$ and
\begin{equation}
        t = \max(\diam C', \diam C'', \dist(C', C'')),
\end{equation}
then
\begin{equation}
        d(x, y) \le t
\end{equation}
for every $x \in C'$ and $y \in C''$, and $C' \cup C''$ is contained
in the closed ball with radius $t$ centered at any point in $C' \cup C''$.
If $C'$, $C''$ are distinct maximal proper sub-cells of a cell $C$,
then $\diam C = t$, because $t \le \diam C$ by the inclusion $C', C''
\subseteq C$, and $t < \diam C$ would imply that there is a proper
sub-cell of $C$ that contains $C'$ and $C''$.  Thus one can take
$\gamma = 1$ when $d(x, y)$ is an ultrametric on $X$ and $\mathcal{C}$
is the cellular structure associated to the ultrametric.

        As in the product case, if $d(x, y)$ satisfies the regularity
conditions (\ref{regularity condition, 1}) and (\ref{regularity
condition, 2}), then $X$ is doubling with respect to $d(x, y)$ if and
only if $X$ is doubling with respect to the cellular structure
$\mathcal{C}$, and a positive Borel measure on $X$ is doubling with
respect to $d(x, y)$ if and only if it is doubling with respect to
$\mathcal{C}$.  If $\widetilde{d}(x, y)$ is another metric on $X$ that
determines the same topology and satisfies the regularity conditions,
then $d(x, y)$ and $\widetilde{d}(x, y)$ are quasisymmetrically
equivalent in the sense that the identity mapping is quasisymmetric as
a mapping from $(X, d(x, y))$ to $(X, \widetilde{d}(x, y))$.

        Of course, there are interesting situations where the
regularity conditions do not hold.  For example, one can have fat
Cantor sets in the real line which the standard Euclidean metric
satisfies (\ref{regularity condition, 1}) and not (\ref{regularity
condition, 2}), and which are doubling with respect to both the metric
and cellular structure.  One may have an upper bound as in
(\ref{regularity condition, 1}) with $\beta < 1$ and not a lower
bound.  It may be that (\ref{regularity condition, 2}) still holds, or
one might ask for a lower bound in terms of a multiple of the
diameters of $C'$ and $C''$.  It may be that the regularity conditions
are satisfied, and that $X$ is quite large and not doubling.

        Let $(X, \mathcal{C})$ be a cellular space, and let $\rho$ be
a nonnegative real-valued function on $\mathcal{C}$ such that $\rho(C)
= 0$ if and only if $C$ has only one element, and
\begin{equation}
        \rho(C') < \rho(C)
\end{equation}
when $C, C' \in \mathcal{C}$, $C' \subseteq C$, and $C' \ne C$.  For
$x, y \in X$, put $d_\rho(x, y) = 0$ when $x = y$, and otherwise
\begin{equation}
\label{d_rho(x, y) = rho(C(x, y))}
        d_\rho(x, y) = \rho(C(x, y))
\end{equation}
where $C$ is the minimal cell that contains $x$ and $y$.  Thus
$d_\rho(x, y) > 0$ when $x \ne y$, and
\begin{equation}
        d_\rho(y, x) = d_\rho(x, y)
\end{equation}
for every $x, y \in X$.  Let us check that
\begin{equation}
        d_\rho(x, z) \le \max(d_\rho(x, y), d_\rho(y, z))
\end{equation}
for every $x, y, z \in X$.  This is trivial when $x = y$ or $y = z$,
and so we may suppose that $x \ne y \ne z$.  If $C(x, y)$, $C(y, z)$
are the minimal cells such that contain $x, y$ and $y, z$,
respectively, then either $C(x, y) \subseteq C(y, z)$ or $C(y, z)
\subseteq C(x, y)$, since $C(x, y)$ and $C(y, z)$ both contain $y$ are
are therefore not disjoint.  Thus $z \in C(x, y)$ or $x \in C(y, z)$,
and the inequality follows.

        This shows that $d_\rho(x, y)$ is an ultrametric on $X$.  The
topology determined by $d_\rho(x, y)$ is the same as the initial
topology on $X$ if for each $x \in X$ and $\epsilon > 0$ there is a
cell $C$ such that $x \in C$ and $\rho(C) < \epsilon$.  By
compactness, this is the same as saying that for each $\epsilon > 0$
there are finitely many cells $C_1, \ldots, C_n$ such that $X =
\bigcup_{i = 1}^n C_i$ and $\rho(C_i) < \epsilon$ for each $i$.

        Suppose that $C$ is a cell with at least two elements.  Thus
$C$ contains a proper sub-cell, and hence a maximal proper sub-cell
$C'$. If $x \in C'$ and $y \in C \backslash C'$, then $C$ is the
minimal cell that contains both $x$ and $y$.  This implies that the
diameter of $C$ is equal to $\rho(C)$ with respect to $d_\rho$, which
holds trivially when $C$ has only one element.  One can check that
each cell $C$ is equal to the closed ball centered at any element of
$C$ with radius $\rho(C)$ with respect to $d_\rho$, and that every
closed ball of positive radius with respect to $d_\rho$ is a cell.

        If $X$ has at least two elements and no isolated points, then
every cell has at least two elements.  This implies that every cell
has at least two distinct maximal proper sub-cells.  In this case, one
can choose $\rho$ so that $d_\rho$ satisfies the regularity conditions
(\ref{regularity condition, 1}) and (\ref{regularity condition, 2}).
As in the earlier examples, one might also be interested in metrics
that do not satisfy the regularity conditions.

        The setting of cellular spaces seems to be quite natural for
having some nice properties while at the same time accommodating
a range of possibilities.

\end{document}